\documentclass{birkjour}

\usepackage{amsthm,amsmath,amssymb,stmaryrd,mathrsfs}
\usepackage{mathtools,amsfonts}
\usepackage{enumerate}
\usepackage{enumitem}
\usepackage{subcaption}

\usepackage{graphicx} 

\usepackage{graphicx}
\usepackage[dvipsnames]{xcolor}
\usepackage{multirow}
\usepackage[hidelinks]{hyperref}
\usepackage{enumerate}
\usepackage{comment}
\usepackage{array}
\usepackage{xspace}
\usepackage{url}

\usepackage{tikz}
\usetikzlibrary{positioning}
\usetikzlibrary{arrows}
\usetikzlibrary{shapes.geometric}

\newcommand\blfootnote[1]{%
  \begingroup
  \renewcommand\thefootnote{}\footnote{#1}%
  \addtocounter{footnote}{-1}%
  \endgroup
}

\newtheorem{theorem}{Theorem}[subsection]
\theoremstyle{definition}
\newtheorem{definition}[theorem]{Definition}
\newtheorem{lemma}[theorem]{Lemma}

\newtheorem{example}[theorem]{Example}

\newcommand{\K}{\ensuremath{{\sf{K}}}\xspace}
\newcommand{\bem}{\ensuremath{{\sf{bem}}}\xspace}
\newcommand{\BEM}{\ensuremath{{\sf{BEM}}}\xspace}
\newcommand{\BM}{\ensuremath{{\mathcal{BM}}}\xspace}
\newcommand{\iKbem}{\ensuremath{{\sf{iK}+\bem}}\xspace}
\newcommand{\iK}{\ensuremath{{\sf{iK}}}\xspace}
\newcommand{\IPC}{\ensuremath{{\sf{IPC}}}\xspace}
\newcommand{\LC}{\ensuremath{{\sf{LC}}}\xspace}
\newcommand{\LP}{\ensuremath{{\sf{LP}}}\xspace}
\newcommand{\BD}{\ensuremath{{\sf{BD}}}\xspace}
\newcommand{\CPC}{\ensuremath{{\sf{CPC}}}\xspace}
\newcommand{\MM}{\ensuremath{{\mathcal{MM}}}\xspace}
\newcommand{\CMM}{\ensuremath{{\mathcal{CMM}}}\xspace}
\newcommand{\BMbem}{\ensuremath{{\mathcal{BM}+\BEM}}\xspace}

\newcommand{\Copy}[2]{{\mathsf{Copy}}(#1,#2)}

\newcommand{\BoxForm}{\ensuremath{{\sf{Form}}_\Box}\xspace}

\newcommand{\prop}{\ensuremath{{\sf{Prop}}\xspace}}
\newcommand{\powerset}[1]{\mathfrak{P}({#1})}
\newcommand{\model}[1]{\mathcal{#1}}
\newcommand{\logic}[1]{\mathcal{#1}}

\title{On the Contingency of Logic in Possible World Semantics}
\author{Iris van der Giessen}
\address{Science Park 107\\
Postbus 94242\\
1090 GE Amsterdam}
\email{i.vandergiessen@uva.nl}

\author{Joost J.\, Joosten, Paul Mayaux, Vicent Navarro Arroyo}
\address{Departament de Filosofia\\
Universitat de Barcelona,\\
Montalegre, 6\\
08001 Barcelona, Catalonia, Spain}
\email{\texttt{jjoosten@ub.edu} \quad \texttt{paul.mayaux@gmail.com} \quad \texttt{vicent.navarro@ub.edu}
}

\date{October 2025}
\subjclass{Primary 03B45; Secondary 03B05, 03B20}
\keywords{Combining logics, Possible worlds semantics, Mixed models, Birelational models}

\begin{document}
\blfootnote{The present paper is the winner of the Spanish Prize of Logic 2025, part of the third edition of the World Logic Prizes Contest.}
\begin{abstract}
    This paper investigates the contingency of logic within the framework of possible world semantics. Possible world semantics captures the meaning of  necessitation, i.e., a statement is necessarily true if it holds in all possible worlds. Standard Kripkean semantics assumes that all possible worlds are governed by one single logic. We relax this assumption and introduce mixed models, in which different worlds may obey different logical systems. The paper provides a first case study where we mix classical propositional logic ($\CPC$) and intuitionistic propositional logic ($\IPC$) in the possible world semantics. We define the class of mixed models $\model M\model M(\CPC, \IPC)$, together with a subclass of concrete mixed models ($\CMM$), and establish their semantic properties. Our main result shows that the set of formulas valid in $\model M\model M(\CPC, \IPC)$ corresponds exactly to the intuitionistic modal logic $\iK$ extended with the Box Excluded Middle axiom ($\iKbem$). To demonstrate this, we prove soundness and completeness results linking $\model M\model M(\CPC, \IPC)$ and $\CMM$, and birelational models for $\iKbem$. 
\end{abstract}
\maketitle
\section{Introduction}
In this paper we will be only concerned with propositional modal logic. Our language thus has letters $p, q, \ldots$ to refer to propositions and we will consider a modality $\Box$ to denote some attribute that such a proposition $p$ can have as in $\Box p$. Typically, the $\Box$ will stand for \textit{necessary} and $\Box p$ should for example be read as `the proposition $p$ is necessarily true'. Modal logic is a system of reasoning that concerns propositional logic together with the modal attribute. The reading conventions for expressions in modal logic follow those of propositional logic where the $\Box$ connective/operator binds as strong as negation.

\subsection{Motivation of our study}
\label{sec:motivation}

Other classical examples to interpret the $\Box$ modality are given by attributes such as \textit{known}, \textit{obligatory} or even more formal notions as \textit{formally provable} giving rise to \textit{epistemic}, \textit{deontic} and \textit{provability logic} respectively \cite{Hintikka1962-HINKAB,Ronnedal2010-RNNAIT,Boolos1993-BOOTLO-7}.

Truth tables are adequate semantics for classical propositional logic. Modal logics require more refined semantics. Ever since the introduction of modal logic, the predominant semantics for modal logic has been so-called \textit{possible world semantics} as introduced by Saul Kripke \cite{Kripke1959-KRIACT, Kripke1963-KRISAO, Kripke1963-KRISCO}. Conceptually the picture of possible world semantics is very clear and suggestive. Certain facts $p$ are contingent and are only true in virtue of some accidental arrangement in the actual world; yet, one can conceive that in other possible worlds $p$ is not true. The necessary facts $q$ are those facts $q$ that cannot be different in any possible world so that $q$ will hold in \textit{any} possible world. Kripke gave the framework of possible world semantics its current flexibility by observing that which worlds are possible may depend on the current world itself. Thus, Kripke introduced a binary accessibility relation $R$ where $aRb$ would stand for $b$ is a possible world of $a$ which we often read as $b$ is \textit{accessible} from $a$. 

One often uses the notation $a\Vdash \varphi$ to denote that the proposition $\varphi$ is \textit{true} (\textit{holds} or more neutral \textit{is forced}) at some world $a$. The forcing condition for statements involving necessity then is defined as 
\[
a\Vdash \Box \varphi \ \ :\Longleftrightarrow \ \ \forall b \ \big( aRb \ \Rightarrow b\Vdash \varphi\big).
\]
In the setting of classical logic, Kripke isolated the system $\K$ of modal logic and proved that it exactly corresponds to the possible world semantics described above. The axioms in this system $\K$ are all classical logical tautologies in the extended language (like $\Box p \to \Box p$) together with all distribution axioms which are axioms of the form $\Box (\varphi \to \psi) \to ( \Box \varphi \to \Box \psi )$.
The reasoning rules in this setting are Modus Ponens $\frac{\varphi \ \ \ \varphi\to \psi}{\psi}$ and Necessitation: $\frac{\varphi}{\Box \varphi}$.

A semantical justification of the rule of Necessitation is readily given. If~$\K$ is to generate those modal formulas that hold in \textit{all} possible worlds then, in particular, all those formulas also hold in \textit{all accessible} worlds. A justification for the rule of Necessitation that does not refer to the possible world semantics could be called an \textit{intensional} or \textit{epistemic justification} and could run as follows. If through pure modal reasoning we arrive at some conclusion $\varphi$ then this reasoning is universal and compelling  whence $\varphi$ is necessary, that is, $\Box \varphi$. The intensional justification seems to presuppose the necessity of reasoning and this assumption is reflected in the possible world semantics in that every possible world obeys the same rules of logic. In this paper we investigate what happens if we drop this assumption. Thus, different possible worlds may obey to different laws of reasoning, to different logics. 

\subsection{General research question}

The object of our study is thus to see what happens with possible world semantics if different logics may govern the reasoning in different possible worlds. We have chosen to associate a possible world $w$ with the set $T_w$ of modal formulas that hold in that world where each $T_w$ will be a theory in a logic $\logic{L}_w$ that may differ from world to world. The modal formulas should still hold by virtue of possible world semantics exclusively.  We thus naturally arrive at the following definition.

\begin{definition}[Mixed Models]
    Let $\{\logic{L}^i\}_{i\in I}$ be a class of logics for some index class $I$. A \emph{mixed model} over $\{\logic{L}^i\}_{i\in I}$ is a triple $\langle W, R, \{T_w\}_{w\in W}\rangle$. Here, $\langle W, R \rangle$ is a modal Kripke frame, that is, $W$ is a non-empty set of possible worlds and $R$ is a binary relation on $W$. 
    
    Each $T_w$ will be a modal theory closed under logical consequence and thus implicitly bears the information of the underlying logic. We will write $\logic{L}_w$ for the logic that corresponds to the world $w$ and we require that $\logic{L}_w \in \{\logic{L}^i\}_{i\in I}$. Furthermore, we shall write $\vdash_w$ instead of $\vdash_{\logic{L}_w}$ for the consequence relation corresponding to the world $w$. The triple $\langle W, R, \{T_w\}_{w\in W}\rangle$ should meet the following requirements.
    \begin{enumerate}
        \item 
        $\bot \notin T_w$;\\
        In other words, each theory $T_w$ is consistent.
        \item 
        $T_w \vdash_w \varphi \ \ \Longrightarrow \ \ \varphi \in T_w$;\\ In other words, each theory $T_w$ is closed\footnote{See Definitions \ref{definition:IPC} and \ref{definition:CPC} for a clarification for the notation $\vdash_w$ in our context.} under $\logic{L}_w$-logical consequence;
        \item 
        $\Box \varphi \in T_w \ \ \Longleftrightarrow \ \ \forall\,  v{\in}W\ (wRv \ \Rightarrow \ \varphi \in T_v)$;\\
        In other words, modal formulas are evaluated through the classic possible world semantics.
    \end{enumerate}
The class of all Mixed Models over $\{\logic{L}^i\}_{i\in I}$ is denoted by $\MM (\{\logic{L}^i\}_{i\in I})$. 
\end{definition}

Given some $\model{M}\in \MM (\{\logic{L}^i\}_{i\in I})$ we often write $\model M , w\Vdash \varphi$ or simply $w\Vdash \varphi$ instead of $\varphi \in T_w$. With $\model M$ being $\langle W, R , \{T_w\}_{w\in W}\rangle$ we will sometimes write $x\in \model M$ instead of $x\in W$. The notation $\model M \models \varphi$ will be short for $\forall \, w {\in}W\ \model M , w\Vdash \varphi$ and $\models \varphi$ will stand for $\forall \, \model M {\in} \MM (\{\logic{L}^i\}_{i\in I}) \ \model M \models \varphi$ if $\{\logic{L}^i\}_{i\in I}$ is clear from the context.

The general research question is to determine a syntactic characterization of $\{ \varphi \mid \forall \, \model M {\in} \MM (\{\logic{L}^i\}_{i\in I}) \ \model M \models \varphi\}$ for a given collection $\{\logic{L}^i\}_{i\in I}$.

\subsection{Narrowing the scope}

This paper will approach the general research question with a first case study. We will here explicitly describe our scope delimitation. First of all, we will consider just two logics: classical propositional logic -- denoted by \CPC-- on the one hand, and intuitionistic propositional logic -- denoted by \IPC-- on the other hand. 

Secondly, we will not directly study $\{\varphi \mid \forall \, \model M{\in} \MM(\{\CPC,\IPC\})\ \model M \models \varphi\}$. Instead, we will look at a smaller class of models $\MM(\CPC,\IPC)$ that are as $\MM(\{\CPC,\IPC\})$ but where we additionally require some behavior on negated $\Box$-formulas.

\begin{definition}\label{definition:ourMixModels}
    The class of models $\MM(\CPC,\IPC)$ are those models $\model M \in \MM(\{\CPC,\IPC\})$ for which we additionally have 
    \begin{enumerate}
    \setcounter{enumi}{3}
        \item 
        $\neg \Box \varphi \in T_w \ \ \Longleftrightarrow \ \ \exists \, v{\in } W \ ( wRv \wedge \varphi \notin T_v).$
    \end{enumerate}
\end{definition}
Let us briefly elaborate on this choice. The idea is that in our first investigation the modal operator is fully governed by possible world semantics in a setting of meta-logic that is classical. We shall show the finite model property which somehow could be seen as circumstantial evidence for this assumption. We thus wish to be both $\Box \varphi$ and $\neg \Box \varphi$ to hold true\footnote{One could consider stipulating $\neg \Box \varphi \in T_w \ \ \Longleftrightarrow \ \ \exists \, v{\in } W \ ( wRv \wedge \neg \varphi \in T_v)$ which would however, yield to a drastically different landscape.} exclusively in virtue of the possible world semantics.  We emphasise that even though the notations are similar with the only difference being the absence or presence of $\{$ and $\}$, the notations $\MM(\{\CPC,\IPC\})$ and $\MM(\CPC,\IPC)$ denote rather different model classes.

\subsection{Plan of the paper}

As agreed upon above, the remainder of this paper shall be dedicated to the study of $\MM(\CPC, \IPC)$ and the corresponding logic. From the outset, it is however not clear that there are any models at all in $\MM(\CPC, \IPC)$ since Clauses 2 -- 4 require more and more formulas to be included in the sets $T_w$ which a priori could make $\bot$ become derivable whence a member of $T_w$ which is excluded by Clause 1. 

To get examples of models, we shall consider the so-called Concrete Mixed Models. These are models where each $T_w$ shall be generated by a semantic model: either a valuation in the case of classical logic or a Kripke model in the case of intuitionistic logic. We will write $\CMM \models \varphi$ to denote that $\varphi$ holds true (in the sense of Definition \ref{definition:ourMixModels}) on all concrete mixed models. 

We will see that the set of those formulas that are true in all models in $\MM(\CPC, \IPC)$ can be associated by an intuitionistic modal logic. Intuitionistic modal logics are defined on an intuitionistic propositional base and usually feature an intuitionistic behavior of $\Box$ \cite{Simpson94PhD}. From Definition \ref{definition:ourMixModels} it should be clear that in our study the $\Box$-formulas will behave classically and thus the law of excluded middle applies to them: $\Box \varphi \vee \neg \Box \varphi$. We define this to be the principle of \textit{Box Excluded Middle} and shall write \bem for it. In Section~\ref{sec:iK}, we will introduce an axiomatization of the intuitionistic modal logic $\iKbem$ together with its semantics in terms of birelational models $\BM$ characteristic for intuitionistic modal logics. We will prove that $\iKbem$ exactly describes the set of all those formulas that are true in all models in $\MM(\CPC, \IPC)$.

In summary, we will prove the following implications as indicated in Figure~\ref{fig:roadmap}. In the conclusion we reflect on our results and provide related literature.

\tikzset{dbl/.style={double,
                     double equal sign distance,
                     implies-implies,
                     shorten >=0pt,
                     shorten <=0pt}}

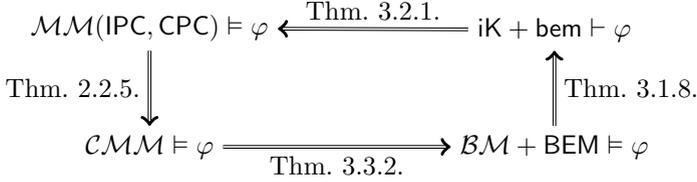
\begin{figure}[ht]
\centering
\begin{tikzpicture}
\node[]        (1)       [] {$\mathcal{MM}(\mathsf{IPC,CPC})\vDash\varphi$};
\node[]        (2)       [below =of 1] {$\mathcal{CMM}\vDash\varphi$};
\node[]        (3)       [right= 3cm of 2] {$\mathcal{BM}+\mathsf{BEM}\vDash\varphi$};
\node[]        (4)       [above=of 3] {$\iKbem\vdash\varphi$};

\draw[->,double] (1) -- (2) node[midway,left] {Thm. 2.2.5.};
\draw[->,double] (2) -- (3) node[midway,below] {Thm. 3.3.2.};
\draw[->,double] (3) -- (4) node[midway,right] {Thm. 3.1.8.};
\draw[->,double] (4) -- (1) node[midway, above] {Thm. 3.2.1.};
\end{tikzpicture}
\caption{Roadmap of the paper.}
\label{fig:roadmap}
\end{figure}

\section{Concrete Mixed Models}

Concrete models will roughly speaking consist of a Kripke model where each node $w$ will either be a model for classical logic or for intuitionistic logic. Before we formally describe concrete mixed models, we shall first settle upon some notation and conventions. 

\subsection{A short recap of intuitionistic logic}

The language that we consider in this paper is the language of modal logic where the modal logical formulas only contain a $\Box$ modality and no $\Diamond$ modality. This is essential since over intuitionistic modal logic these are non-dual. Thus, the formulas that we consider will be denoted $\BoxForm$ and are generated by
\[
\mathcal F \ \ := \ \ \top \mid \bot \mid  \prop \mid (\mathcal F \wedge \mathcal F) \mid (\mathcal F \vee \mathcal F) \mid (\mathcal F \to \mathcal F) \mid \Box \mathcal F .
\]
We use the notation $\neg \varphi$ as shorthand for $(\varphi \to \bot)$. We inherit bracketing conventions from propositional logic where $\Box$ now binds equally strong as $\neg$. Propositional logical formulas are generated similarly by omitting the $\Box \mathcal F$ case. Let us recall a well-known Hilbert-style presentation of Intuitionistic Propositional Calculus.

\begin{definition}[\IPC]\label{definition:IPC}
The axioms of \IPC are all formulas of one of the forms:
\begin{align*}
    &A \rightarrow (B \rightarrow A) && A \rightarrow \big(B \rightarrow (A \land B)\big)\\
    &(A \rightarrow B) \rightarrow \Big(\big(A \rightarrow (B \rightarrow C)\big) \rightarrow (A \rightarrow C)\Big) && (A \land B) \rightarrow A\\
   &A \rightarrow (A \lor B) && (A \land B) \rightarrow B \\
   &B \rightarrow (A \lor B) &&\bot \rightarrow A \\
   &(A \lor B) \rightarrow \Big((A \rightarrow C) \rightarrow \big((B \rightarrow C) \rightarrow C\big)\Big)  && 
\end{align*}
The sole rule of inference is Modus Ponens, allowing one to conclude $B$ from proofs of $A$ and $A\to B$ where proofs are as usual finite sequences of formulas each being either an axiom or arising as the conclusion from an application of the Modus Ponens rule applied to two earlier elements of the sequence.
\end{definition}

We will use the notation of $\IPC$ both for the calculus formulated in \BoxForm and in the purely propositional language; the context will always clarify which scenario applies. We will write $\Gamma \vdash_\IPC \varphi$ to indicate that there exists an \IPC-proof of $\varphi$ from $\Gamma$. Namely, that there is some finite sequence each element of which is either an axiom of \IPC, a member of $\Gamma$, or the conclusion of an application of Modus Ponens to two formulas that occur earlier in the sequence. Additionally, the last member of the sequence should be $\varphi$.  Again, the context should specify where the formulas in $\Gamma, \varphi$ are supposed to range over.

\begin{definition}[\CPC]\label{definition:CPC}
Classical logic will be denoted by \CPC and arises from \IPC by adding the single axiom scheme of the excluded middle: $\neg A \lor A$.
\end{definition}

Again, the context will specify where the formulas in the schemes range over and $\Gamma \vdash_\CPC \varphi$ is defined similarly as $\Gamma \vdash_\IPC \varphi$.
Let us now recall Kripke semantics for intuitionistic logic.
\begin{definition}
An \emph{intuitionistic Kripke model} is a tuple $ \model M := \langle W, \leq, V\rangle$ where $W$ is a non-empty set with $\leq$ being a reflexive, transitive, and antisymmetric relation on $W$ and $V : W \to \powerset \prop$ is a valuation on $W$ that is required to be monotone in the sense that $x\leq y \Rightarrow V(x) \subseteq V(y)$.
\end{definition}

For an intuitionistic Kripke model $\model M := \langle W, \leq, V\rangle$, by reading convention, we will write $x \in \model M$ instead of $x\in W$ and likewise for other tuples. We say that $\model M$ is \emph{rooted} whenever there is a root $r\in W$ so that $r\leq x$ for any $x\in W$. We use the same notation $\Vdash$ for forcing in intuitionistic Kripke models that is inductively defined in the standard way: 
\begin{align*}
    x\Vdash p &\ \ :\Leftrightarrow \ \ p{\in}V(x);\\
    x\nVdash \bot;\\
    x\Vdash \varphi \vee \psi&\ \ :\Leftrightarrow \ \ x\Vdash \varphi \mbox{ or }x\Vdash \psi;\\
    x\Vdash \varphi \wedge \psi&\ \ :\Leftrightarrow \ \ x\Vdash \varphi \mbox{ and }x\Vdash \psi;\\
    x\Vdash \varphi \to \psi&\ \ :\Leftrightarrow \ \ \forall \, y{\geq}x\, (\mbox{ if } y\Vdash \varphi \mbox{, then }y\Vdash \psi).
\end{align*}

We state the well-known soundness and completeness for \CPC and \IPC and refer to \cite{vanDalen1980-VANLAS-2} for details.

\begin{theorem}
For $\varphi$ in the language of propositional logic we have
\begin{align*}
    \vdash_\IPC \varphi \ \ &\Longleftrightarrow \ \ \mbox{each rooted Kripke model $\model M$ with root $r$ satisfies } \model M, r \Vdash \varphi.\\
    \vdash_\CPC \varphi \ \ &\Longleftrightarrow \ \ \mbox{for any one-point Kripke model $\model M$ we have } \model M \models \varphi.
\end{align*}

\end{theorem}

\subsection{Concrete Mixed Models}

We can now define the set $\CMM$ of \emph{Concrete Mixed Models}.

\begin{definition}
A Concrete Mixed Model is a triple $\langle W, R, \{ \model M_w\}_{w\in W} \rangle$ where $\langle W,R \rangle$ is a Kripke frame ($W$ not empty and $R\subseteq W\times W$) and each $\model M_w = \langle W_w,\leq_w,V_w\rangle$ is a rooted intuitionistic Kripke model. Also, all the $\model M_w$ will have disjoint domains. The root of $\model M_w$ will always be denoted by $\overline w$ and the class of all concrete mixed models by  $\CMM$. 
\end{definition}

By definition, the domains of all the $\model M_w$'s are disjoint. Thus, we can safely write $\leq$ instead of $\leq_w$ and likewise for $V$ instead of $V_w$. For a concrete mixed model $\model M$ we will write $x\in \model M$ to indicate that $x$ is an element of some of the $\model M_w$'s. It is clear that for each $x\in \model M$ there is a unique $w\in W$ so that $w\leq x$. We define $\model M, x \Vdash \varphi$ by recursion as follows\footnote{By now $\Vdash$ stands for forcing in modal models, in mixed models, in intuitionistic Kripke models and in mixed models alike. The context will always determine the exact notion.}.

\begin{definition} 
Let $\model M \in \CMM$ and $x\in \model M$. We define the forcing relation $\Vdash$ as usual for the atomic propositions and propositional connectives: 
\begin{align*}
    x\Vdash p &\ \ :\Leftrightarrow \ \ p{\in}V(x);\\
    x\nVdash \bot;\\
    x\Vdash \varphi \vee \psi&\ \ :\Leftrightarrow \ \ x\Vdash \varphi \mbox{ or }x\Vdash \psi;\\
    x\Vdash \varphi \wedge \psi&\ \ :\Leftrightarrow \ \ x\Vdash \varphi \mbox{ and }x\Vdash \psi;\\
    x\Vdash \varphi \to \psi&\ \ :\Leftrightarrow \ \ \forall \, y{\geq}x\, (\mbox{ if } y\Vdash \varphi \mbox{, then }y\Vdash \psi).
\end{align*}
Finally, for the modality we stipulate: 
\[
x\Vdash \Box \varphi :\Leftrightarrow \Big( \forall v \, \big( w R v \Rightarrow \overline v \Vdash \varphi \big) \mbox{ and $w$ is the unique world so that $\overline w \leq x$}\Big).
\]
\end{definition}

The intuition behind the concrete mixed models is that the root $\overline w$ of each Kripke model $\model M_w$ will generate a theory $T_w$ of those formulas that are forced at the root. The collection of those theories $T_w$ then defines a mixed model. Before we prove this we make some basic observations on concrete mixed models. Firstly, we observe that the truth of $\Box$-formulas is quite well-behaved.

\begin{lemma}\label{theorem:BasicCMMLemma}
Let $\model M \in \CMM$ and $x,y \in \model M$ and let $\varphi \in \BoxForm$. We have:
\begin{enumerate}
\item
$x\leq y \ \ \Longrightarrow \ \ \big( x\Vdash \Box \varphi \Leftrightarrow y\Vdash \Box \varphi \big)$;
\item
$x\Vdash \varphi \ \& \ x\leq y \ \ \Longrightarrow \ \ y\Vdash \varphi$;
\item
$\vdash_\IPC \varphi \ \ \Longrightarrow \ \ x\Vdash \varphi$;
\item
If $\model M_x$ consists of the single world $\overline x$, then $\vdash_\CPC \varphi \ \ \Longrightarrow \ \ x\Vdash \varphi$. 
\end{enumerate}
\end{lemma}

\begin{proof}
The first item is clear since both $x$ and $y$ share the same unique $\overline w $ so that $\overline w \leq x$ and $\overline w \leq y$. Thus, as far as intuitionistic semantics is concerned, $\Box$-formulas behave as propositional variables that either hold on the entire model or do not hold on the entire model. Consequently, the second item is proved by an easy induction on the complexity of $\varphi$. The third item follows by a simple induction on the length of an \IPC proof of $\varphi$. Given that $\Box$-formulas are monotonous (preserved under $\leq$) the proof is as for propositional ones in \IPC. Via a simple induction on $\varphi \in \BoxForm$ we see that $\overline w \Vdash \varphi \vee \neg \varphi$  when $\model M_w$ consists of the single world $\overline w$ that establishes the fourth item. 
\end{proof}

We now see that concrete mixed models provide a recipe to obtain mixed models.

\begin{theorem}\label{thm:MM-existence}
Let $\model M := \langle W, R, \{ \model M_w\}_{w\in W} \rangle$ be a concrete mixed model. For each $w\in W$, we define $T_w : = \{ \varphi \in \BoxForm \mid \overline w \Vdash \varphi \}$. The triple $\langle W, R, \{ T_w\}_{w\in W} \rangle$ now defines a mixed model of $\MM (\CPC, \IPC)$.
\end{theorem}

\begin{proof}
If $\model M_w$ consists of the single world $\overline w$, we set $\logic L_w := \CPC$. However, $T_w$ can define a classical theory also in other ways, so, in general, we set $\logic L_w := \CPC$ whenever $p\vee \neg p \in T_w$ for all $p\in \prop$ 
and we set  $\logic L_w := \IPC$ otherwise. We need to check the four items from Definition \ref{definition:ourMixModels}. The first item is trivial. The second item follows from the previous lemma. The third item follows from the definition of the~$T_w$ and the fourth item follows from Lemma \ref{theorem:BasicCMMLemma}.
\end{proof}

We observe that for concrete mixed models we have that $\varphi \vee \psi \in T_w$ implies that either $\varphi \in T_w$ or $\psi \in T_w$. Thus, it is not hard to see that there are mixed models that are not induced by a concrete mixed model.

Also observe that the constructed mixed model in Theorem~\ref{thm:MM-existence} forces exactly the same formulas as the considered concrete mixed model. So, concrete mixed models are just a special case of mixed models and thus we have the following theorem. This proves the left-most arrow of Figure~\ref{fig:roadmap}.

\begin{theorem}[Soundness with respect to Concrete Mixed Models]
\quad

    \[
      \MM (\CPC, \IPC) \models \varphi \ \ \ \Longrightarrow \ \ \CMM  \models \varphi.
    \]
\end{theorem}

Furthermore, we observe that our move from $\MM (\CPC, \IPC)$ to $\CMM$ suggests a more general strategy for $\MM (\{\logic L^i\}_{i\in I})$: replace each world $w$ in the modal Kripke frame by a structure of semantics that is sound and complete for $\logic L_w$ so that this structure generates a set closed under $\logic L_w$ derivability.

\section{A logic for Mixed Models}

To axiomatize the set of formulas true in all models of $\MM (\CPC, \IPC)$ we look at intuitionistic modal logics. The reason to investigate intuitionistic modal logic in our case study is that models of $\CPC$ can be considered as one-world intuitionistic Kripke models of $\IPC$, as discussed in the previous section.

\subsection{Intuitionistic modal logic}
\label{sec:iK}

The general idea of intuitionistic modal logic is to capture the meaning of the modalities in an intuitionistic manner. There are different motivations and strategies to do so, e.g., see \cite{BozicDosen84,FischerServi77,Simpson94PhD}. In our context, we work with the $\Box$-modality only. One way to define an intuitionistic counterpart of classical modal logic~$\K$ (see Section \ref{sec:motivation}) is as follows:

\begin{definition}
    Intuitionistic modal logic $\iK$ is defined to be the minimal set of formulas containing the intuitionistic axioms from Definition~\ref{definition:IPC} together with the modal axiom 
    \[
    \Box (\varphi \to \psi) \to ( \Box \varphi \to \Box \psi ),
    \]
    and is closed under rule Modus Ponens $\frac{\varphi \ \ \ \varphi\to \psi}{\psi}$ and Necessitation: $\frac{\varphi}{\Box \varphi}$.
\end{definition}

A standard way to define meaning for intuitionistic modal logics is via birelational semantics. This features two relations: a relation~$\leq$ to model intuitionistic reasoning and a relation $R$ to model modal reasoning.

\begin{definition}
\label{def:birelational}
A \emph{birelational frame} is a tuple $ \model F := \langle W, \leq, R\rangle$ so that  
$\langle W, \leq\rangle$ is an intuitionistic Kripke frame and so that 
    $x\leq y \ \wedge \ y R z \ \Longrightarrow \ xRz$. 
    
    A \emph{birelational model} is a tuple $ \model M := \langle W, \leq,  R, V\rangle$ so that  
$\langle W, \leq,  R\rangle$ is a birelational frame and $V : W \to \powerset \prop$ is a valuation on $W$ that is required to be monotone in the sense that $xRy \Rightarrow V(x) \subseteq V(y)$.
    The class of birelational models is denoted $\BM$.

    The forcing relation on birelational models is defined as follows:
    \begin{align*}
    x\Vdash p &\ \ :\Leftrightarrow \ \ p{\in}V(x);\\
    x\nVdash \bot;\\
    x\Vdash \varphi \vee \psi&\ \ :\Leftrightarrow \ \ x\Vdash \varphi \mbox{ or }x\Vdash \psi;\\
    x\Vdash \varphi \wedge \psi&\ \ :\Leftrightarrow \ \ x\Vdash \varphi \mbox{ and }x\Vdash \psi;\\
    x\Vdash \varphi \to \psi&\ \ :\Leftrightarrow \ \ \forall \, y{\geq}x\, \big(\mbox{ if } y\Vdash \varphi \mbox{, then }y\Vdash \psi\big);\\
    x\Vdash \Box \varphi &\ \ :\Leftrightarrow \ \ \forall y \big(\mbox{ if } xRy\mbox{, then }y\Vdash \varphi\big).\\
\end{align*}
\end{definition}

\begin{figure}[t!]
    \centering
    \begin{subfigure}[t]{0.4\textwidth}
        \centering
        \includegraphics[height=1.2in]{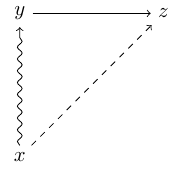}
        \caption{The frame condition for the monotonicity of the valuation.}
        \label{fig:subfig:monotonicity_fc}
    \end{subfigure}%
    ~ 
    \begin{subfigure}[t]{0.4\textwidth}
        \centering
        \includegraphics[height=1.2in]{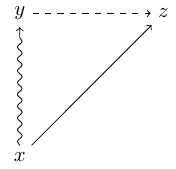}
        \caption{The \BEM frame condition.}
        \label{fig:subfig:BEM_fc}
    \end{subfigure}
    \caption{The straight arrows represent the relation $R$ and the snake arrow the relation $\leq$. The dashed arrow indicates the relations we should have, given the undashed arrows.}
\end{figure}

The frame condition $x\leq y \ \wedge \ y R z \ \Longrightarrow \ xRz$ makes sure that the models are monotone in $\leq$, meaning that $x\leq y \ \wedge \ x \Vdash \phi \ \Longrightarrow \ y \Vdash \phi$. In Figure \ref{fig:subfig:monotonicity_fc} we depict this frame condition.
    It is well-known (\cite{BozicDosen84}) that $\BM$ is sound and complete for $\iK$. 
\begin{theorem}\label{iks}
    $\iK$ is sound and complete with respect to $\BM$.
\end{theorem}

For our case study, we already observed that $\Box$-formulas behave classically, namely, they satisfy the \emph{Box Excluded Middle} axiom $\Box \varphi \vee \neg \Box \varphi$, which we call $\bem$. For this reason, we will investigate intuitionistic modal logic with this additional axiom.

\begin{definition}
    Logic $\iKbem$ is the minimal set of formulas containing $\iK$ and $\bem$ and is closed under Modus Ponens and Necessitation.
\end{definition}

To describe this logic semantically, we introduce a new frame property which we will call $\BEM$. 

\begin{definition}
    The acronym \BEM~(Box Excluded Middle) stands for the frame condition 
    \[
    \forall x,y,z\ \big( x \leq y \ \wedge \ xRz \to y Rz\big ). 
    \]
    The class of birelational models satisfying $\BEM$ is denoted by $\BM + \BEM$.
\end{definition}
In Figure \ref{fig:subfig:BEM_fc} we represent the \BEM frame condition. It is easy to see that \BEM is a sufficient condition for \bem to hold.
\begin{lemma}
    Let $\model F$ be a birelational frame $ \model F := \langle W, \leq, R\rangle$. We have that
    \[
    \model F \mbox{ satisfies } \BEM \ \ \Longrightarrow \ \ \model F \models \bem .
    \]
\end{lemma}

\begin{figure}[t]
        \centering
        \includegraphics[width=0.3\linewidth]{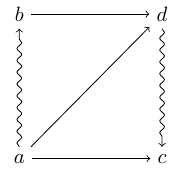}
        \caption{Counterexample. The straight arrows represent the  relation $R$ and the snake arrows the relation $\leq$.}
        \label{fig:counter}
    \end{figure}

\begin{proof}
With $x \in \model F$ arbitrary we assume $x\nVdash \Box A$ and set out to prove that $x\Vdash \neg \Box A$. So, we consider $y$ arbitrary with $x\leq y$ and wish to see that $y\nVdash \Box A$. But by the assumption that $x\nVdash \Box A$, there is $z$ with $xRz$ and $z\nVdash A$. By the assumption that $x \leq y \ \wedge \ xRz \to y Rz$ we conclude that $yRz$ so that indeed $y\nVdash \Box A$. Finally, since $y$ was arbitrary, we can conclude that $x\Vdash \neg\Box A$.
\end{proof}
Indeed, this lemma shows that \BEM is sufficient for \bem to hold. This is the only direction that we need. The converse implication in this lemma is not true as the following counterexample demonstrates. 
\begin{example}
    Let us consider a frame with four points $a,b,c$ and $d$. The points $a,b,c$ are related as $aRc$ and $a\leq b$ and without $bRc$. Moreover, we put $bRd$, $aRd$ and $d\leq c$. Figure \ref{fig:counter} shows the non-trivial relations in our model of four points. We observe that this frame does not satisfy $\BEM$. However, it can be shown that for each valuation on this frame we have $\Box \varphi\vee \neg\Box \varphi$ forced in every world of the frame, for any $\varphi$.

    Firstly, we observe that we only need to prove the claim for $a$ and $b$. On the one hand, assume $d\Vdash \varphi$, for some arbitrary $\varphi$. By monotonicity, we have that $c\Vdash \varphi$. Therefore, $a\Vdash\Box \varphi$ and $b\Vdash\Box\varphi$, so $a\Vdash \Box\varphi\vee\neg\Box\varphi$ and $b\Vdash \Box\varphi\vee\neg\Box\varphi$.

    On the other hand, assume $d\nVdash \varphi$. Then, $a\nVdash\Box\varphi$ and $b\nVdash\Box\varphi$. By the definition of the intuitionistic negation, we derive that $a\Vdash\neg \Box\varphi$ and $b\Vdash\neg\Box\varphi$, so $a\Vdash\Box \varphi\vee\neg\Box\varphi$ and $b\Vdash \Box\varphi\vee\neg\Box\varphi$. 
\end{example}

Thus, we see that \BEM is not a necessary condition for \bem to hold. The following theorem is proved by a standard canonical model construction. 

\begin{theorem}[Intuitionistic modal completeness]\label{theorem:iKbemComplete}\ 
    \[
  \BM + \BEM \models \varphi  \ \ \Longrightarrow \ \ \iKbem \vdash \varphi.
    \]
\end{theorem}
\begin{proof}
     To define the canonical model we consider \emph{prime} sets of formulas. A set of formulas $\Gamma$ is \emph{prime} if it satisfies the following conditions:
    \begin{itemize}
        \item (Deductively closed): $\Gamma \vdash_{\iKbem} \varphi \Longrightarrow \varphi \in \Gamma$;
        \item (Consistent): $\Gamma \not \vdash_{\iKbem} \bot$;
        \item (Disjunction property): $\varphi \vee \psi \in \Gamma \Longrightarrow \varphi \in \Gamma \mbox{ or } \psi \in \Gamma$. 
    \end{itemize}
    We define the following birelation model $\model M = \langle W,\leq,R,V\rangle$ with
    \begin{align*}
        W \ &= \ \ \{ \Gamma \mid \Gamma \mbox{ is prime} \};\\
        \Gamma \leq \Gamma' \ &:\Leftrightarrow  \ \Gamma \subseteq \Gamma';\\
        \Gamma R \Delta \ &:\Leftrightarrow  \ \ \mbox{if } \Box \varphi \in \Gamma, \mbox{ then } \varphi \in \Delta;\\
        V(\Gamma) \ &\coloneqq \ \ \{ p \mid p \in \Gamma \}.
    \end{align*}
    Model $\model M$ satisfies frame property $\BEM$. To see this suppose $\Gamma \leq \Gamma'$ and $\Gamma R \Delta$. Suppose that $\Box \varphi \in \Gamma'$. Since $\Gamma$ is deductively closed over $\iKbem$, we have $\Box \varphi \vee \neg \Box \varphi \in \Gamma$. By the disjunction property we have $\Box \varphi \in \Gamma$ or $\neg \Box \varphi \in \Gamma$. By the facts that $\Gamma \subseteq \Gamma'$ and that $\Gamma'$ is consistent we have $\Box \varphi \in \Gamma$. Now, by definition of $\Gamma R\Delta$ this means that $\varphi \in \Delta$, as desired. Therefore, $\model M \in \BMbem$. 

    By an easy induction on $\varphi$, we can prove the following truth lemma:
    \[
    \Gamma \models \varphi \Longleftrightarrow \varphi \in \Gamma.
    \]

    Now, to prove the theorem, suppose $\iKbem \not \vdash \varphi$. 
    By a standard Lindenbaum type argument we can see that there is a prime set $\Gamma$ such that $\varphi \notin \Gamma$. By the truth lemma, we conclude that $\model M, \Gamma \not \Vdash \varphi$. Therefore, $\BMbem \not \models \varphi$.
\end{proof}

So we have established the right arrow of our roadmap in Figure~\ref{fig:roadmap}. It remains to prove the upper and lower arrows in Figure~\ref{fig:roadmap}, performed in the subsequent two sections.

\subsection{From logic $\iKbem$ to Mixed Models}

It is easy to see that mixed models are sound for the logic \iKbem.

\begin{theorem}[Soundness with respect to Mixed Models]
\quad

    \[
    \iKbem \vdash \varphi \ \ \ \Longrightarrow \ \  \MM (\CPC, \IPC) \models \varphi.
    \] 
\end{theorem}
\begin{proof}
    The prove of this soundness result is quite standard. First, checking that all intuitionistic axioms from Definition~\ref{definition:IPC} are valid in all mixed models is straightforward but tedious. By Definition \ref{definition:ourMixModels}, the $\bem$ property is trivial. It is easy to check that Modus Ponens preserves validity. Finally, we develop the reasoning for the preservation of validity under Necessitation. Suppose $\varphi$ is valid in $\MM(\CPC,\IPC)$ and fix any mixed model $\model M$ and world $w\in \model M$. For every $v$ with $wRv$, we have that $\varphi\in T_v$, by the validity of $\varphi$. Then, $\Box\varphi\in T_w$ and $\Box\varphi$ holds at $w$. Because $\model M$ and $w$ were arbitrary, we get the validity of $\Box\varphi$ in $\MM(\CPC,\IPC)$. 
\end{proof}
So we have established the upper arrow in Figure~\ref{fig:roadmap}.

\subsection{From Concrete mixed models to birelational models}

It is quite easy to see that each concrete mixed model gives rise to a birelational one.

\begin{theorem}\label{thm:3.3.1}
Let $\model M \in \CMM (\CPC, \IPC)$ with $\model M = \langle W, R, \{ \model M_w\}_{w\in W} \rangle$ and each $\model M_w = \langle W_w, \leq_w, V_w \rangle$. We define the birelational model $\model M' = \langle W', R', \leq', V' \rangle$ by setting 
\begin{align*}
        W' \ &:= \ \ \bigcup_{w\in W}W_w;\\
        xR'y \ &:\big( \overline w R \overline v \ \wedge \ \overline w \leq_w x \ \wedge \ \overline v = y \big);\\
        x\leq'y \ &:\Leftrightarrow  \ \ x\leq_w y \text{ for some }w\in W;\\
        V(v) \ &\coloneqq \ \ \{ p \mid \overline w \leq_w v \wedge p \in V_w(v) \}.
    \end{align*}
The birelational model $\model M'$ is in $\BMbem$. Moreover, the concrete mixed model~$\model M$ makes true the same modal formulas as the birelational model $\model M'$ in the sense that
\[
\forall \, w{\in}\model M\, \forall \, \varphi{\in}\BoxForm \ \big(  \model M, w \Vdash \varphi \ \Leftrightarrow    \model M', w \Vdash \varphi \big).
\]
\end{theorem}

\begin{proof}
We first observe that $R'$ is well-defined since for each $x\in \bigcup_{w\in W}W_w$ there is exactly one $w\in W$ with $\overline w \leq_w x$. Likewise, $\leq'$ is well-defined since for each pair $x,y \in \bigcup_{w\in W}W_w$ there is at most one $w\in W$ so that $x\leq_w y$. We can check that it satisfies the frame property $x\leq' y \ \wedge \ y R' z \ \Longrightarrow \ xR'z$ and that it is monotone in $\leq'$. Therefore, $\model M'$ is a birelational model. 

We now prove that $\model M'$ satisfies frame property $\BEM$. Suppose $x \leq' y$ and $xR'z$. This means $x \leq_w y$ for some $w \in W$ and, therefore, $\overline w \leq_w x$ and $\overline w \leq_w y$. From $xR'z$ it follows that there is some $\overline v$ such that $\overline w R \overline v$, $\overline w \leq_w x$ and $\overline v = z$. Now, also $\big( \overline w R \overline v \ \wedge \ \overline w \leq_w x \ \wedge \ \overline v = y \big)$ and, thus, $yRz$ as desired.

By induction on $\varphi$ we can now prove that $\model M, w \Vdash \varphi \ \Leftrightarrow    \model M', w \Vdash \varphi$. Since $\leq' \  = \ \leq$, the only interesting case is when $\varphi = \Box \psi$. Here we observe that $xR'y \Leftrightarrow (\overline w\leq x \wedge wRv \wedge \overline  v = y)$ is exactly what is needed for the inductive step. 
\end{proof}
Observe that the model $\model M'$ makes $\bem$ true because it satisfies the frame property $\BEM$, by Theorem 3.2.1 and the fact that model $\model M'$ makes the same theorems true as model $\model M$.

To conclude our proof and close Diagram \ref{fig:roadmap} we  need the more complicated construction that translates an arbitrary birelational model into a concrete mixed model. This construction is presented in the proof of the following theorem. 


\begin{theorem}\label{theorem:FromCMMToiKbem}\ 
    \[
  \CMM \models \varphi \ \ \Longrightarrow \ \ \BM + \BEM \models \varphi  
    \]
\end{theorem}

\begin{proof}
    By contraposition, let us fix a birelational model 
    $\model M = \langle W, R, \leq , V\rangle$ so that, at some $x\in \model M$, we have $\model M, x \nVdash \varphi$ for some particular formula $\varphi$. We shall see how to come up with a concrete mixed model $\model M'$ and some $x' \in \model M'$ so that $\model M' , x'\nVdash \varphi$. It is good to recall that the $\Vdash$ for mixed models  means something quite different from $\Vdash$ in relational semantics.

In the construction, we consider for each world $w\in W$ a copy for all its $\leq$-visible worlds. To make this precise, for each $x,y \in W$ such that $x \leq y$, we introduce element $\Copy{x}{y}$. We use these elements to define a concrete mixed model
\[
\model M' = \langle W' , R', \{ \model N_x\}_{x\in W'}\rangle
\]
where the constituent parts will be defined now.
We define
\[
W' := \{ 
\Copy{x}{x} \mid x\in W
\}
\]
and 
\[
\Copy{x}{x} R' \Copy{y}{y} \ \ :\Longleftrightarrow \ \ xRy .
\]
Note that $W$ and $W'$ have the same size.
Finally, we define for each $\Copy{x}{x} \in W'$ the intuitionistic Kripke model \[
\model N'_{\Copy{x}{x}} := \langle W'_{\Copy{x}{x}},\ \leq'_{\Copy{x}{x}},\ V'_{\Copy{x}{x}} \rangle
\]
where 
\[
W'_{\Copy{x}{x}} \ := \ \{ \Copy{x}{y} \mid x,y \in W \text{ and } x\leq y\}
\]
where the root of $W'_{\Copy{x}{x}}$ is $\overline{\Copy{x}{x}} := \Copy{x}{x}$.
We next define for $\Copy{x}{y}, \Copy{x}{y'} \in W'_{\Copy{x}{x}}$
\[
\Copy{x}{y} \leq'_{\Copy{x}{x}} \Copy{x}{y'} \ \ :\Longleftrightarrow \ \ y\leq y'.
\]

Note that all $\leq'_{\Copy{x}{x}}$ are disjoint so we may just as well write $\leq'$ as we will often do.
Finally, we define $V'_{\Copy{x}{x}}$ by stipulating that for all $\Copy{x}{y} \in W'_{\Copy{x}{x}}$,
\[
p \in V'_{\Copy{x}{x}} (\Copy{x}{y}) \ \ :\Longleftrightarrow \ \ p\in  V(y).
\]
Indeed, it is not hard to see that each  $\model N'_{\Copy{x}{x}}$ defines an intuitionistic model. The reflexivity, transitivity and antisymmetry of $\leq'_{\Copy{x}{x}}$ are given by the fact that $\leq$ is a poset because $\model M$ is a birelational model. Finally, to check the monotonicity of the valuation $V'_{\Copy{x}{x}}$, we assume
$\Copy{x}{y}\leq'_{\Copy{x}{x}}\Copy{x}{y'}$ and $p\in V'_{\Copy{x}{x}}(\Copy{x}{y})$. Then, we try to prove that $p\in V'_{\Copy{x}{x}}(\Copy{x}{y'})$.  By hypothesis, we deduce that $y\leq y'$ and $p\in V(y)$. By the monotonicity of $\leq$, we derive that $p\in V(y')$ and, by the definition of $V'_{\Copy{x}{x}}$, we conclude that $p\in V'_{\Copy{x}{x}}(\Copy{x}{y'})$.

By induction on $\psi$, one can now prove
\[
\model M , x \Vdash \psi \ \ \Leftrightarrow \ \ \model M', \Copy{x}{x} \Vdash \psi .
\]
For the induction to work, the  proof will actually prove that 
\[
\model M , x \Vdash \psi \ \ \Leftrightarrow \ \ \model M', \Copy{a}{x} \Vdash \psi,
\]
whenever $a \leq x$ in $\model M$. The only non-trivial case is the $\Box$ modality. We should prove that
\[\mathcal{M},x\Vdash \Box\varphi \ \ \Leftrightarrow \ \ \mathcal{M'},\Copy{a}{x}\Vdash \Box\varphi. \]
For the left to right implication, assume that $\mathcal{M},x\Vdash\Box\varphi$ and consider $\Copy{y}{y}\in W'$ such that $\Copy{a}{a}R'\Copy{y}{y}$, where $\Copy{a}{a}$ is the unique world so that $\overline{\Copy{a}{a}}\leq'_{\Copy{a}{a}}\Copy{a}{x}$. We ought to prove that $\model M',\overline{\Copy{y}{y}}\Vdash\varphi$. By the assumption that $\model M, x \Vdash \Box \varphi$, we know that for all elements $z \in W$ such that $xRz$ we have $\model M, z \Vdash \varphi$. Now, since $\Copy{a}{a}R'\Copy{y}{y}$ we have, by construction of $R'$, that $aRy$. From $a \leq x$ and $aRy$ it follows, by $\BEM$, that $xRy$, so $\model M, y \Vdash \varphi$. By induction hypothesis, we know that $\model M',\Copy{y}{y} \Vdash \varphi$, as desired. 


For the right to left implication, assume that $\mathcal{M'},\Copy{a}{x}\Vdash\Box\varphi$ with $a \leq x$ and consider some $y\in W$ such that $xRy$. Our aim is to prove that $\model M, y\Vdash\varphi$. From $a \leq x$ and $xRy$ it follows that $a R y$ (because $\model M$ is an intuitionistic birelational model (Definition~\ref{def:birelational})). This means that $\Copy{a}{a}R'\Copy{y}{y}$. Now, by the assumption that $\mathcal{M'},\Copy{a}{x}\Vdash\Box\varphi$, we know that $\model M',\Copy{y}{y} \Vdash \varphi$ since $\overline{\Copy{a}{a}}\leq'\Copy{y}{y}$. By inductive hypothesis, we conclude that $\model M,y\Vdash\varphi$, as we wanted to prove. 
\end{proof}

This finishes the final direction of our roadmap in Figure~\ref{fig:roadmap}.

\section{Conclusions and future work}

This paper departs from a philosophical debate by asking the question how one could define necessity in a possible world semantics when dealing with different logical reasoning systems assigned to each world. To the best of our knowledge, this paper presents one of the first treatments of mixed models in its current form. Departing from a similar philosophical debate, \cite{freire_martins_2024} introduces models in which the logical reasoning in each world is described by lattice structures. Closely related to the concrete mixed models introduced in the present paper is the work in \cite{MaderiaMartinsBenevides2019} in which each world of an epistemic (S5)-frame is assigned a model from a so-called knowledge representation framework. Strictly speaking, \cite{MaderiaMartinsBenevides2019} does not mix different logics, but that could similarly be said for our case study in which a CPC-model could be considered as a one-point IPC-model. Also related is the work in \cite{DiaconescuStefaneas2007,Martins_etal2011} in a categorical setting. A recent study that is rather similar in spirit to ours but departs from different motivations has been presented in 
    \cite{mojtahedi2022provabilitylogicha} and in \cite{mojtahedi2025provabilitymodels}. 
    
 Substantial research to other forms of mixing logics has been done, for which we only mention the following few: fusions of modal logics \cite{Gabbay_etall03}, fibring logics \cite{Gabbay_1996,Gabbay1998}, combining classical and constructive connectives/modalities in different ways such as in ecumenical logic \cite{Binder2022, Williamson1988-WILEAE, delCerro1996,  marin2020ecumenicalmodallogic}.

In this paper, we introduced a novel framework to study what happens when  different logics govern the reasoning in different possible worlds. We developed the notion of mixed models and by focusing on $\CPC$ and $\IPC$, we gave concrete instantiations of mixed models. We established that the valid formulas in this setting are precisely captured by the intuitionistic modal logic $\iKbem$. This result highlights that the coexistence of distinct logics within a unified possible world semantics can be rigorously modeled, while still admitting a clean syntactic characterization. It seems that \iKbem has the finite model property: Theorem \ref{theorem:iKbemComplete} seems to yield finite models by restricting only to the relevant formula and its subformulas. In turn, Theorem \ref{theorem:FromCMMToiKbem} preserves finiteness. However, the Devil is in the details and things would need thorough checking.

Future work could extend this investigation in several directions. One natural step is to study mixed models involving larger families of logics, such as intermediate logics or substructural logics. An interesting question arises when we consider mix incomparable logics, by which we mean to mix logics  that are not contained in each other. A concrete example would be to study the case when we consider the Gödel-Dummett logic --\LC axiomatized over \IPC via the scheme $(p\to q)\vee (q\to p)$-- and the intuitionistic logic of depth bounded by 2, $\BD_2$ axiomatized by $p\vee \big(p\to (q\vee \neg q)\big)$.

Furthermore, it would be interesting to study the same research question when we consider the family of many-valued logics and, in particular, three-valued systems like the strong Kleene logic, $\K_3$, or the Logic of Paradox, \LP. Another direction is to investigate properties on the modal relation between said worlds such as transitivity, reflexivity, etc. It also makes sense to combine the set-up in the realm of \textit{temporal logics} since time may behave differently globally (governed by the $R$ relation) than locally (governed by the $\leq$ relation).

Finally, the connections between mixed models and philosophical debates on logical pluralism and the necessity of reasoning deserve further examination. By broadening the scope of possible world semantics in this way, we aim to contribute to a deeper understanding of the contingency of logic itself.
\section*{Acknowledgments}
We thank the anonymous referees for their valuable comments to improve the paper. We thank Lev Beklemishev, David Fernández Duque, Johan van Benthem, Genoveva Martí, José Martínez Fernández and Sven Rosenkranz for interesting discussions. Joosten acknowledges the support of the Generalitat de Catalunya (Government of Catalonia)
under grants 2022 DI 051 (Departament d’Empresa i Coneixement) and 2021 SGR
00348 and support of the Spanish Ministry of Science and Innovation under grant
PID2023-149556NB-I00, Dynamics of Gödel Incompleteness (DoGI). Van der Giessen acknowledges the support by the Dutch Research Council (NWO) under the project \textit{Finding Interpolants: Proofs in Action} with file number VI.Veni.232.369 of the research programme Veni.
\bibliographystyle{abbrv}
\bibliography{journal_version}
\end{document}